\documentclass{article}
\usepackage{amsthm}
\usepackage{amsfonts}

\begin{document}

Tsemo Aristide

College Boreal

1 Yonge Street Toronto, ON

M5E 1E5

tsemo58@yahoo.ca

\bigskip
\bigskip

\centerline{\bf  Linear Foliations on affine manifolds.}

\bigskip
\bigskip

\centerline{\bf Abstract.}

{\it In this paper, we study  affine manifolds endowed with linear foliations. These are foliations defined by vector subspaces invariant by the linear holonomy. We show that an $n$-dimensional  compact, complete, and oriented affine manifold endowed with a codimension $1$ linear foliation ${\cal F}$ is homeomophic to the $n$-dimensional torus if the  leaves of ${\cal F}$ are simply connected. Let $(M,\nabla_M)$ be  a $3$-dimensional compact affine manifold endowed with a codimension $1$ linear foliation. We prove that $(M,\nabla_M)$ has a finite cover which is homeomorphic to the total space of a bundle over the circle if its developing map is injective, and has a  convex image. }

\bigskip


\bigskip

{\bf A.M.S. Classification codes.}

53C12, 53C15.
\bigskip
\bigskip

\centerline{\bf 1. Introduction.}

An affine manifold $(M,\nabla_M)$ is a differentiable manifold $M$ endowed with a flat, torsion-free connection $\nabla_M$, or equivalently with an atlas whose coordinate changes are affine transformations. 
An affine manifold is unimodular if it has a parallel volume form.  In 1962, Markus [21] has conjectured that a compact, complete, and oriented affine manifold $(M,\nabla_M)$ is geodesically complete if and only if it is unimodular. Remark that the Markus conjecture implies that a compact affine manifold is either aspherical, or its first Betti number is nonzero (see Goldman and Hirsch [15] p. 644).  The Markus conjecture has been proven by:
 Carriere [6] if the holonomy of $(M,\nabla_M)$ preserves a flat Lorentzian metric.
  Goldman and Hirsch have shown that the Markus conjecture is true if the holonomy of $(M,\nabla_M)$ is solvable and its solvable rank is inferior to the dimension of $M$,  or if $(M,\nabla_M)$ is endowed with a rational Riemannian metric (see [16] p.176, p.191).
  The existence of parallel forms on a compact affine manifold is studied  by Goldman and Hirsch in [15]. They have shown that if the holonomy of $(M,\nabla_M)$ is nilpotent, and $M$ has a parallel $k$-form whose cohomology class is non zero, for every $j\leq k$ there exists a non zero parallel $k$-form (see [15] p.645). The existence of a nonzero parallel $1$-form is shown in [13] (theorem 4.1) if $(M,\nabla)$ is compact not radiant, and its holonomy group is nilpotent. The kernel of a non zero parallel $k$-form integrates to a foliation of codimension $k$ which is an example of linear foliation studied in this paper. Linear foliations are examples of analytic foliations studied by Haefliger [17], [18] and transversely affine foliations studied by Bobo Seke [5]. 
 
 In [31],  I have shown that the restriction of the connection $\nabla_M$ to $aff(M,\nabla_M)$, the Lie algebra of $Aff(M,\nabla_M)$, the group of affine transformations of  $(M,\nabla_M)$ defines an associative structure. This enables to define non trivial linear foliations if $Aff(M,\nabla_M)$ is not discrete.  
 
In this paper, we define affine manifolds and their group of automorphisms  in section 2. 
Let $(M,\nabla_M)$ be a closed $3$-dimensional manifold whose group of automorphisms $Aff(M,\nabla_M)$ is not discrete. The associative algebra $aff(M,\nabla_M)$  is the direct sum of a nilpotent associative algebra $N_M$, and a semi-simple associative algebra $S_M$. In section 3, the proposition 3.2  shows that if $N_M$ is non zero and the holonomy representation is injective, either there exists a translation which commutes with the holonomy of $M$, or    a finite cover of $M$ is the total space of a bundle over the circle. The proposition 3.3 shows that if $S_M$ is non zero, either $(M,\nabla_M)$ is radiant, or $(M,\nabla_M)$ is endowed with a linear foliation of codimension $1$. The proofs of these results use the associative structure of $aff(M,\nabla_M)$ and  the classification of $3$-dimensional manifolds whose fundamental group is solvable.
 
Let $(M,\nabla_M)$ be an $n$-dimensional compact and oriented affine manifold endowed with a linear foliation ${\cal F}_M^U$ of codimension $1$. In section $4$, the corollary 4.1 shows that if $(M,\nabla_M)$ is complete and ${\cal F}_M^U$   has a leaf whose fundamental group is solvable, then $(M,\nabla_M)$ is unimodular. The proposition 4.4 shows that if the leaves of ${\cal F}_M^U$ are simply connected and $(M,\nabla_M)$ is complete, then $M$ is homeomorphic to the $n$-dimensional torus. The proof of this result is obtained  by remarking that the fundamental group of $M$ is isomorphic to the global holonomy of the foliation which is commutative, and by using the rigidity theorem of Malcev (see [25]).  The propositions 4.5 and 4.8 show that if $M$ has a compact leaf without holonomy, or if $(M,\nabla_M)$ is complete and has a compact leaf, then a finite cover of $M$ is the total space of a bundle over $S^1$. These results are obtained  by applying the fact that the closure of a union of compact leaves is a union of compact leaves (Haefliger [18]), and using  the fact  that the holonomy of a compact and complete affine manifold is irreducible shown by Fried Goldman and Hirsch in [13]. We  remark  that if $(M,\nabla_M)$ is compact and radiant, the Bobo Seke cohomology class of ${\cal F}_M^U$ does not vanish.  

 In section 5, the corollary 5.1 shows the following result: let $(M,\nabla_M)$ be a compact $3$-dimensional affine manifold whose developing map is injective, and has a convex image. Suppose that $(M,\nabla_M)$ is endowed with codimension $1$ linear foliation. Then a finite cover of  $M$ is the total space of a bundle over $S^1$. To obtain this result, we prove that the fundamental group of $M$ is solvable by using results of Stallings on group theory, Novikov on foliations, and we use the classification of $3$-dimensional manifolds whose fundamental group is solvable.

\bigskip

\centerline{\bf Plan.}

\bigskip

1. Introduction.

2.  Affine manifolds and their developing maps.

3. Linear foliations: first properties.

4.  Codimension 1 linear foliations.

5. Codimension $1$ linear foliations and $3$-dimensional affine manifolds.

\bigskip

{\bf 2. Affine manifolds and their developing maps.}

\medskip

All manifolds in this paper will be assumed to be connected and without boundary unless otherwise stated.

\medskip

{\bf Definition 2.1. The developing map.}
Let $(M,\nabla_M)$ be an affine manifold, and $\hat M$ its universal cover, the pullback $\nabla_{\hat M}$ of $\nabla_M$ to $\hat M$ defines on $\hat M$ the structure of an affine manifold. There exists a local diffeomorphism $D_M:\hat M\rightarrow \mathbb{R}^n$, such that a chart containing an element $ x'\in \hat M$ can be defined by an open subset $\hat U$ of $\hat M$ which contains $ x'$, and such that the restriction ${D_M}_{\mid\hat U}$ of $D_M$ to $\hat U$ is a diffeomorphism onto its image. The map $D_M$ is called the developing map.

\medskip

{\bf Definitions 2.2.}
Let $(M,\nabla_M)$  (resp. $(M',\nabla_{M'}))$ an $n$-dimensional affine manifold (resp. an $n'$-dimensional affine manifold) defined by the affine atlas $(U_i,\phi_i)_{i\in I}$ (resp. $(V_j,\psi_j)_{j\in J})$. A morphism of affine manifolds $f:(M,\nabla_M)\rightarrow (M',\nabla_{M'})$ is a differentiable map $f$, such that $\psi_j\circ f_{\mid U_i}\circ \phi_i^{-1}:\mathbb{R}^n\rightarrow\mathbb{R}^{n'}$ is an affine map.
We denote by $Aff(M,\nabla_M)$  the group of affine transformations of $(M,\nabla)$. It is the group of  diffeomorphisms of $M$ which preserve $\nabla_M$. We denote  by $Aff(M,\nabla_M)_0$ the connected component of the identity of $Aff(M,\nabla_M)$.

\medskip

Let $Gl(\mathbb{R}^n)$ (resp. $Aff(\mathbb{R}^n)$)  be the group of linear transformations (resp. affine transformations) of the vector space $\mathbb{R}^n$. For every element $g\in Aff(\mathbb{R}^n)$, there exists an unique $L(g)\in Gl(\mathbb{R}^n)$, and $c_g\in\mathbb{R}^n$, such that for every element $x$ of $\mathbb{R}^n$,  $g(x)=L(g)(x)+c_g$. The map $L(g)$ is called the linear part of $g$, and $c_g$ the translational part of $g$. We will identify $g$ with $(L(g),c_g)$.
 The developing map defines a representation $H_M:Aff(\hat M,\nabla_{\hat M})\rightarrow Aff(\mathbb{R}^n)$ such that the following diagram is commutative:

$$
\matrix{ \hat {M}&{\buildrel{g}\over{\longrightarrow}}& \hat M\cr D_M\downarrow &&\downarrow D_M\cr M&{\buildrel{H_M(g)}\over{\longrightarrow}}& M}
$$

The restriction $h_M$ of $H_M$ to $\pi_1(M)$ is called the holonomy. For every element $\gamma$ of $\pi_1(M)$, we can write $h_M(\gamma)=(L(h_M)(\gamma),a_{\gamma})$. The representation $L(h_M)$ is called the linear holonomy of $(M,\nabla_M)$. It is the holonomy of the connection $\nabla_M$ in the classical sense. 

The group $Aff(M,\nabla_M)$ is the quotient of the normalizer $N(\pi_1(M))$ of $\pi_1(M))$ in $Aff(\hat M,\nabla_{\hat M})$ by $\pi_1(M)$. The Lie algebra $aff(M,\nabla_M)$  of $Aff(M,\nabla_M)$ is isomorphic to the Lie algebra $n(\pi_1(M))$ of $N(\pi_1(M))$.
 In  [31], I have shown that the restriction of $\nabla_M$ to $aff(M,\nabla_M)$  is an associative product, and the  morphism of Lie algebras $Lie(H_M):aff(\hat M,\nabla_{\hat M})\rightarrow aff(\mathbb{R}^n)$ associated to $H_M$ is a morphism of associative algebras where $aff(\mathbb{R}^n)$ is endowed with the canonical structure of associative algebra defined by $(A,a).(B,b)=(AB,A(b))$.

\medskip

{\bf Definitions 2.3.}
An affine manifold $(M,\nabla_M)$ is complete if and only if the connection $\nabla_M$ is geodesically complete. This is equivalent to saying that the developing map is a diffeomorphism. 
The affine manifold $(M,\nabla_M)$ is radiant if $h_M$ fixes an element of $\mathbb{R}^n$. This is equivalent to saying that there exists a translation $t_a$ such that $t_ah_Mt_{-a}=L(h_M)$.
The affine manifold is unimodular if and only if the image of $L(h_M)$ is contained in $Sl(n,\mathbb{R})$ the group of linear automorphisms whose determinant is $1$.

\medskip

{\bf Definitions 2.4.}
A $p$-dimensional  affine submanifold $(F,\nabla_F)$ of the affine manifold $(M,\nabla_M)$ is a $p$-dimensional submanifold $F$ of $M$ endowed with the structure of an affine manifold, such that the canonical embedding $i_F:(F,\nabla_F)\rightarrow (M,\nabla_M)$ is a morphism of affine manifolds.

\medskip

Let $\hat F$ be the universal cover of $F$, we can lift $i_F$ to an affine map $\hat i_F:\hat F\rightarrow \hat M$. The image of $D_M\circ \hat i_F$ is contained in a $p$-dimensional affine subspace $U_F$ of $\mathbb{R}^n$. The map $D_M\circ \hat i_F:\hat F\rightarrow U_F$ is a developing map of $F$.
There exists a canonical  morphism $\pi_F:\pi_1(F)\rightarrow \pi_1(M)$ induced by $i_F$. Let $\gamma$ be an element of $\pi_1(F)$, the holonomy $h_F(\gamma)$ is the restriction of $h_M(\pi_F(\gamma))$ to $U_F$. If there is no confusion, we are going to denote $h_M(\pi_F(\gamma))$ by $h_M(\gamma)$.
Remark that the kernel of $h_F$ can be non trivial and the kernel of $h_M$ trivial as shows the following example:

\medskip

{\bf Example 2.1.}
Consider the quotient of $\mathbb{R}^n-\{0\}$ by $h(x)=2x$, $n>2$, it is an $n$-dimensional affine manifold $(M,\nabla_M)$, its universal cover is $\mathbb{R}^n-\{0\}$. Let $V$ be a $2$-dimensional vector subspace of $\mathbb{R}^n$. The quotient of $V-\{0\}$ by $h$ is a $2$-dimensional affine submanifold of $(M,\nabla_M)$ diffeomorphic to the $2$-dimensional torus $T^2$. Its holonomy has a non trivial kernel.

\bigskip

{\bf 3. Linear foliations: first properties.}

\medskip

The purpose of this section is to define and to study some general properties   of linear foliations.

\medskip

{\bf Definitions 3.1.}
Let $(M,\nabla_M)$ be a $n$-dimensional affine manifold,  $h_M:\pi_1(M)\rightarrow Aff(\mathbb{R}^n)$ its holonomy representation, and $L(h_M)$ its linear part. Suppose that there exists a $p$-dimensional vector subspace $U$ of $\mathbb{R}^n$ such that $L(h_M)(\pi_1(M))(U)=U$. Then, the foliation of $\mathbb{R}^n$ whose leaves are  $p$-dimensional affine subspaces  parallel to $U$ is invariant by $h_M(\pi_1(M))$. The pullback of this foliation by $D_M$ is a foliation ${\cal F}^U_{\hat M}$, defined on $\hat M$ invariant by $\pi_1(M)$, which is the pullback of a foliation  ${\cal F}^U_M$ defined on $M$ by the covering map $p_M:\hat M\rightarrow M$.
We will say that ${\cal F}^U_M$ is a linear foliation.
Since $U$ is invariant by the linear holonomy, $h_M$ passes down to a representaion $h_{{\cal F}_M^U}$ on the quotient which is called the global holonomy of   the linear foliation ${\cal F}_M^U$. 

\medskip

{\bf Examples 3.1.}
Let $A$ be an hyperbolic element of $Gl(2,\mathbb{Z})$, that is, an integral $2\times 2$  matrix which has two real eigenvalues $\lambda, {1\over\lambda}$ such that $|\lambda|>1$. We can define $T_A$ the suspension of the $2$-dimensional torus $T^2$ over the circle $S^1$, which is the quotient of $\mathbb{R}\times \mathbb{R}^2$ by the group generated by the affine transformations $f_1,f_2,f_3$ defined by $f_1(x,y,z)=(x+1,A(y,z))$, $f_2(x,y,z)=(x,y+1,z)$ and $f_3(x,y,z)=(x,y,z+1)$.  Let $(y_0,z_0)$ an eigenvector associated to the eigenvalue $\lambda$. The  plane $P_{\lambda}$ generated by $(1,0,0)$ and $(0,y_0,z_0)$ is invariant by the linear holonomy of $T_A$ and defines a linear foliation ${\cal F}_{T_A}^{P_{\lambda}}$ on $T_A$. The leaves of this foliation can be cylinders, planes or mobius bands. Similarly, one can define the linear foliation ${\cal F}_{T_A}^{P_{1\over\lambda}}$.  Ghys and Sergiescu [14] have shown that every infinitely differentiable codimension $1$  foliation defined on $T_A$ which is $C^1$  close to ${\cal F}_{T_A}^{P_{\lambda}}$ (resp. $C^1$ close to ${\cal F}_{T_A}^{P_{1\over\lambda}}$) is $C^{\infty}$ conjugated to it.

Let $(M,\nabla_M)$ be an affine manifold, $Lie(H_M)(n(\pi_1(M))$, the image of $n(\pi_1(M))$ by $Lie(H_M)$ is stable by the associative product of $aff(\mathbb{R}^n)$. The theorem $23$  III  of Albert [1] implies that  we can write $Lie(H_M)(n(\pi_1(M))=S_M\oplus N_M$ where $S_M$ is a semi-simple associative algebra, and $N_M$ a nilpotent associative algebra.
$S_M=\sum_{i=1}^{i=l}M(n_i,\mathbb{K}_i), \mathbb{K}_i$ is the field of real, complex or quaternion numbers.

\medskip

{\bf Proposition 3.1.}
{\it Let $(M,\nabla_M)$ be an $n$-dimensional  affine manifold such that $Aff(M,\nabla_M)$ is not discrete. Then $(M,\nabla_M)$ is endowed with a non trivial linear foliation or $(M,\nabla_M)$ is a radiant affine manifold.}

\medskip

{\bf Proof.}
The fact that $Aff(M,\nabla_M)$ is not discrete implies that either $S_M$ or $N_M$ is not the zero algebra.

Suppose that $N_M$ is not the zero algebra, there exists a non zero  element $(C,c)\in N_M$ such that $(C,c)^2=(C^2,C(c))=0$. If $C\neq 0$, then the linear holonomy preserves $Ker(C)$, the kernel of $C$ and $Im(C)$, the image of $C$ which are proper, non trivial, vector subspaces of $\mathbb{R}^n$. We deduce the existence of linear foliations ${\cal F}_M^{Ker(C)}$ and ${\cal F}_M^{Im(C)}$ on $(M,\nabla_M)$.
If $C=0$, then $c\neq 0$. For every element $\gamma\in\pi_1(M), L(h_M)(\gamma)(c)=c$. We deduce that the line $D$, generated by $c$ defines a linear foliation ${\cal F}_M^D$ on $(M,\nabla_M)$.

Suppose that $S_M$ is not the zero algebra, there exists a non zero idempotent $E=(E',e')$ in $S_M$. We have $(E',e')^2=({E'}^2,E'(e'))=(E',e')$; this implies that for every $x\in\mathbb{R}^n, E(x)=E'(x+e')$. We deduce that we can change the origin and suppose that $E$ is linear.  The vector subspace $Ker(E)$  is preserved by the linear holonomy. If $E$ is not bijective, then there exists a linear foliation on $(M,\nabla_M)$ defined by ${\cal F}_M^{Ker(E)}$. If $E$ is a bijective idempotent, $E$ is the identity map and its flow are homothetic maps.  The fact that the holonomy commutes with the flow of $E$ implies that it preserves the origin, and $(M,\nabla_M)$ is a radiant affine manifold.

\medskip

{\bf Proposition 3.2.}
{\it Let $(M,\nabla_M)$ be a $3$-dimensional compact affine manifold such that $h_M$ is injective, suppose that $N_M$ is not trivial. Then, either there exists a non trivial translation which commutes with the holonomy, or a finite cover of
 $M$ is the total space of a bundle over $S^1$.}

\medskip 

Firstly, we prove:

\medskip

{\bf Lemma 3.1.}
{\it Let $(M,\nabla_M)$ be a closed  $3$-dimensional affine manifold whose fundamental group is solvable. Then there exists a finite cover of $M$ which is the total space of a bundle over $S^1$.}

\medskip

{\bf Proof.}
The classification of manifolds with a solvable fundamental group (see [2] Theorem 1.20) implies that a finite cover of $M$ is the total space  of a bundle over $S^1$, or $M$ is homeomorphic to $P\mathbb{R}^3\#P\mathbb{R}^3$ which cannot be endowed with an affine structure (see Smillie [28] Theorem 1).

{\bf Proof of proposition 3.2.}
Suppose that $N_M$ is not trivial, there exists a non zero element $(C,c)$ of $N_M$ such that $(C,c)^2=(C^2,C(c))=(0,0)$, if $C=0$. We have seen in the proof of the proposition 3.1 that the translation generated by $c$ commutes with the holonomy.
If $C\neq 0$, then $dim(Im(C))$, the dimension of $Im(C)$ is $1$, $dim(Ker(C))=2$ and $Im(C)\subset Ker(C)$. Since the linear holonomy preserves $Im(C)$ and $Ker(C)$, we deduce that the holonomy $h_M(\pi_1(M))$ is solvable. This implies that $\pi_1(M)$ is solvable since $h_M$ is injective.

\medskip

{\bf Proposition 3.3.}
{\it Let $(M,\nabla_M)$ be a $3$-dimensional compact affine manifold, suppose that $S_M$ is not trivial. Then either $(M,\nabla_M)$  is radiant or there exists a codimension $1$ linear foliation defined on $(M,\nabla_M)$.}

\medskip

{\bf Proof.}
There exists a non zero idempotent $E$ which is an element of $S_M$. We can suppose that $E$ is linear (see the proof of the proposition 3.1).  Suppose that $Ker(E)=0$. Then $E$ is the identity. Since the holonomy $h_M$ commutes with the  flow of $E$, we deduce that $(M,\nabla_M)$ is a radiant manifold.
If $Ker(E)\neq 0$, either $dim(Ker(E)=2$, or $dim(Im(E)=2$, we deduce that either ${\cal F}_M^{Ker(E)}$ or ${\cal F}_M^{Im(E)}$ is a codimension $1$ linear foliation.

\bigskip

{\bf 4. Some properties of codimension $1$ linear foliations.}

\medskip

In this section, we are going to study some properties of codimension $1$ linear foliations.

\medskip

{\bf Proposition 4.1}
{\it Let $(M,\nabla_M)$ be an affine manifold  whose developing map is injective, and which is endowed with a linear foliation  ${\cal F}_M^U$ of codimension $1$. Suppose that ${\cal F}_M^U$ has a leaf $F_0$ whose fundamental group is solvable, and there exists a leaf $\hat F_0$ of ${\cal F}_{\hat M}^U$, above $F_0$ preserved by $Ker(h_{{\cal F}_M^U})$. Then $\pi_1(M)$ is solvable.}

\medskip

{\bf Proof.}
We can identify $\pi_1(M)$ with the image of the holonomy since the developing map is injective. Consider  $r:Ker(h_{{\cal F}_M^U})\rightarrow Aff(\hat F_0)$ the restriction of $Ker(h_{{\cal F}_M^U})$, the kernel of the global holonomy to $\hat F_0$.  Let $(e_1,...,e_n)$ be a basis of $\mathbb{R}^n$ such that $(e_1,...,e_{n-1})$ is a basis of $U$. Let $\gamma$ be an element of the kernel of $r$, $L(\gamma(e_i))=e_i, i<n$, and $L(\gamma)(e_n)=e_n+u_\gamma, u_\gamma\in U$ since $\gamma$ is an element of $Ker(h_{{\cal F}_M^U})$. We deduce that  the kernel of $r$ is an unipotent group.  The image of $r$ is solvable since $\pi_1(F_0)$ is solvable. This implies  that $Ker(h_{{\cal F}_M^U})$ is a solvable group. We deduce that  $\pi_1(M)$ is solvable since $h_{{\cal F}_M^U}(\pi_1(M))\subset Aff(\mathbb{R})$.

\medskip

{\bf Proposition 4.2.}
{\it Let $(M,\nabla_M)$ be an affine manifold whose developing map is injective. Suppose that $(M,\nabla_M)$ is endowed with a linear foliation ${\cal F}_M^U$  of codimension $1$. Let $F_0$ be a leaf of ${\cal F}_M^U$, and $\hat F_0$ be a leaf of ${\cal F}_{\hat M}^U$ above  $F_0$. We suppose that there exists  $x_0\in\mathbb{R}^n$, such that $\hat F_0=(x_0+U)\cap D_M(\hat M)$. Then $\hat F_0$ is preserved by $Ker(h_{{\cal F}_M^U})$.}

\medskip

{\bf Proof.}
Suppose that $x_0+U$ contains $\hat F_0$. In general, $(x_0+U)\cap D_M(\hat M)$ can have several connected components, and $Ker(h_{{\cal F}_M^U})$ acts on the set of these connected components.  Since we have assumed that it is equal to $\hat F_0$,  it is connected and preserved by the group  $Ker(h_{{\cal F}_M^U})$.

 \medskip
 
 In [16] p.176, Fried, Goldman and Hirsch have shown the following result:
 
 \medskip
 
 {\bf Proposition 4.3.}
 {\it Let $(M,\nabla_M)$ be a compact affine manifold with solvable holonomy group $\Gamma$. Then:
 
 a) If $(M,\nabla_M)$ is complete, then $M$ has a parallel volume form.
 
 b) If $(M,\nabla_M)$ has a parallel volume form and the solvable rank of $\Gamma$ is less than the dimension of $M$, then $M$ is complete.}
 
 \medskip
 
 From this result, we deduce:

\medskip

{\bf Corollary 4.1.}
{\it Let $(M,\nabla_M)$ be a compact oriented and complete affine manifold. Suppose that $(M,\nabla_M)$ is  endowed with a codimension $1$ linear foliation ${\cal F}^U_M$ which has a leaf whose fundamental group is  is solvable. Then $(M,\nabla_M)$ is unimodular.}

\medskip

\medskip

{\bf Corollary 4.2.}
{\it  Let $(M,\nabla_M)$ be a  compact and unimodular affine manifold. Suppose that $(M,\nabla_M)$ is   endowed with a codimension $1$ linear foliation ${\cal F}_M^U$ which has a leaf $F_0$ whose fundamental group is solvable. Then $(M,\nabla_M)$ is complete if the developing map is injective,  there exists a leaf $\hat F_0$ of ${\cal F}_{\hat M}^U$ above $F_0$ such that $\hat F_0=(x+U)\cap D_M(\hat M)$,  and the solvable rank of the holonomy is inferior or equal to the dimension of $M$.}

\medskip

{\bf Proof of corollaries 4.2 and 4.3.}
The propositions 4.1 and 4.2 imply that $\pi_1(M)$ is solvable. We conclude by using proposition 4.3.

\medskip

{\bf Proposition 4.4.}
{\it Let $(M,\nabla_M)$ be an $n$-dimensional compact, oriented, and complete affine manifold endowed with a linear foliation ${\cal F}_M^U$ of codimension $1$. Suppose that all the leaves of ${\cal F}_M^U$ are simply connected.  Then $M$ is homeomorphic to the $n$-dimensional torus $T^n$.}

\medskip

{\bf Proof.}
Let $\gamma$ be an element of $\pi_1(M)$ and $p:\mathbb{R}^n\rightarrow \mathbb{R}^n/U$ the quotient map. Suppose that the linear part of $h_{{\cal F}_M^U}(\gamma)$ is not the identity of $\mathbb{R}^n/U=\mathbb{R}$. This implies that if fixes an element $x$ of $\mathbb{R}$. We have $\gamma(p^{-1}(x))=p^{-1}(x)$ and the fundamental group of $p_M(p^{-1}(x))$ is not trivial. Contradiction. We deduce that  the image of $h_{{\cal F}_M^U}$ is commutative. Let $x$ be an element of $\mathbb{R}$, and  $\gamma$  an element of $Ker(h_{{\cal F}_M^U})$, the kernel of $h_{{\cal F}_M^U}$ that we suppose distinct of the identity. The fact that $h_{{\cal F}_M^U}(\gamma)(x)=x$ implies also that $\gamma(p^{-1}(x))=p^{-1}(x)$, since the fundamental group of $p_M(p^{-1}(x))$ is trivial, we deduce that the restriction of $\gamma$ to $p^{-1}(x)$ is the identity. This is in contradiction with the fact the action of $h_M(\gamma)$ on $\mathbb{R}^n$ is free.  We deduce that $\pi_1(M)$ is commutative since it is isomorphic to  the image of $h_{{\cal F}_M^U}$. The theorem in [16] p. 194,  implies that $(M,\nabla_M)$ is a nilmanifold. Since $\pi_1(M)$ is commutative, the corollary 2, p. 34 of [25] implies that $M$ is homeomorphic to the $n$-dimensional torus $T^n$.

\medskip

{\bf Remark 4.1.}
At  page 190 of  [16], Goldman and Hirsch  have shown that the Markus conjecture is true if $Aff(M,\nabla_M)$ acts transitively on $M$. I have generalized this result by showing that the Markus conjecture is true if $dim(Aff(M,\nabla_M)\geq dim(M)-1$.( See Tsemo [31]). The previous proposition is related to Tsemo [31] since in Tsemo [31], it is shown that the orbits of $Aff(M,\nabla_M)_0$ are leaves of a linear foliation of codimension $1$ if $dimAff(M,\nabla_M)_0=dim(M)-1$, and if $(M,\nabla_M)$ is compact and complete or compact and unimodular.

Let $M$ be a closed manifold endowed with a $C^2$ foliation whose leaves are diffeomorphic to $\mathbb{R}^{n-1}$. Rosenberg [26] has shown that the fundamental group of $M$ is free Abelian and the universal cover of $M$ is contractible. In [26], he has shown that if the dimension of $M$ is $3$, then $M$ is homeorphic to the $3$-dimensional torus. If the dimension of $M$ is superior or equal to $4$, results of Wall and Hsiang [20] and in Freedman and Quinn [11] implies that $M$ is homeomorphic to the $n$-dimensional torus.

The Borel conjecture asserts that two aspherical compact manifolds with isomorphic fundamental group are homeomorphic. Farrell and Jones [10] corollary 0.2 have shown that if $(M,\nabla_M)$ and $(N,\nabla_N)$ are compact complete affine manifolds with isomorphic fundamental groups, they are homeomorphic if $n\neq 3,4$.

Let ${\cal F}$ be a codimension $1$ analytic foliation defined on the closed manifold $M$. Haefliger (see [ [18] p. 386), has shown that the number of compact leaves of ${\cal F}$ is finite, or all the leaves of ${\cal F}$ are compact.  We adapt the arguments of Haefliger in the context of affine manifolds:

\medskip

{\bf Proposition 4.5.}
{\it Let $(M,\nabla_M)$ be a compact affine manifold endowed with a codimension $1$ oriented linear foliation ${\cal F}_M^U$. Suppose that $M$ has a compact leaf $F_0$ without holonomy. Then every leaf of ${\cal F}_M^U$ is compact and there exists an affine morphism $(M,\nabla_M)\rightarrow S^1$ whose fibres are the leaves of ${\cal F}_M^U$.}

\medskip

{\bf Proof.}
The Reeb stability theorem implies that $C$, the union of compact leaves without holonomy is open; $C$ is also closed. To show that,  consider $x$ an element of the closure $\bar C$ of $C$. Consider a connected affine chart $V$ which contains $x$ and an affine transversal $T$ to ${\cal F}_M^U\cap V$ which contains $x$ and an element $y$ of $C$; $T$ can be constructed with a geodesic which contains $x$ and $y$.  Let $F_x$ be the leaf of ${\cal F}_M^U$ containing $x$. The first theorem of the paragraph 3.2 of Haefliger [18] p. 386 implies that the leaf $F_x$ is closed. We show now that $F_x$ does not have holonomy. For every element $\gamma$ of $\pi_1(F_x)$, the holonomy of $\gamma$ at $T$ is an affine map $h_\gamma$ which  fixes $x$. We can identify it to the restriction  of an element of $Gl(\mathbb{R})$.  Since ${\cal F}_M^U$ is oriented, the proposition 3.1 of [18] implies that we can suppose that $F_y\cap T=\{y\}$. This implies that $h_\gamma(y)=y$. We deduce that $h_\gamma$ is the identity since it  is the restriction of a linear map of $\mathbb{R}$ which fixes two points.  This implies that the holonomy group $F_x$ is trivial,  and $C$ is open and closed. We deduce  that $C=M$.

The space of leaves of ${\cal F}^U_M$ is endowed with an affine structure whose charts are defined by local affine transversal.

\medskip

We  study now the leaves which have an infinite holonomy:

\medskip

{\bf Proposition 4.6.}
{\it Let $(M,\nabla_M)$ be a compact affine manifold endowed with a codimension $1$  linear foliation ${\cal F}_M^U$, which has  compact leaves $F_1,...,F_p$, such that the holonomy of $F_1$ is infinite. We suppose that:

(i). The developing map is a covering, thus $(M,\nabla_M)$ is the quotient of $D_M(\hat M)$ by $h_M(\pi_1(M))$;

(ii). there exists leaves $F'_1,...,F'_p$ of ${\cal F}_{\hat M}^U$,  such that $p_M(F'_i)=F_i$ and $\cup_{i=1,..,p}D_M(F'_i)=U\cap D_M(\hat M)$.

Then, there exists a basis $(e_1,...,e_n)$ of $\mathbb{R}^n$ such that $(e_1,...,e_{n-1})$ is a basis of $U$ and for every $\gamma\in\pi_1(M)$, $h_M(\gamma)(0)\in U$. }

\medskip

{\bf Proof.}
 Let $\gamma$ be an element of $\pi_1(F_1)$ which has an infinite holonomy. Write $L(h_M))(\gamma)(e_n)=ae_n+u, u\in U$, we can suppose that $|a|<1$ and up to a change of origin, that $h_M(\gamma)(0)\in U$.
Let $\gamma'\in\pi_1(M)$, suppose that $h_M(\gamma')(0)=be_n+v, v\in U$, for every integer $l$, $h_M(\gamma^l\gamma')(0)=a^lbe_n+u_l, u_l\in U$.  We are going to show that $b=0$. Suppose that $b\neq 0$.

Let $x$ be an element of $F_1$ and $x'\in  F'_1$ such that $p_M(x')=x$. Consider  open subsets $V$ of $M$ and $V'$ of $\hat M$ such that $x\in V, x'\in V'$, the restriction of $p_M$ to $V'$ is a diffeomorphism onto $V$, and the restriction of $D_M$ to $V'$ is a diffeomorphism onto its image. Since $lim_la^lb=0$, there exist an integer $l$ and an element $z"\in U$ such that $h_M(\gamma^l\gamma')(z")\in D_M(V')$. Since $U\cap D_M(\hat M)=\cup_{i=1,..,p}D_M( F'_i)$, we deduce that there exists $i\in\{1,..,p\}, z'\in F'_i$ such that $D_M(z')=z"$, write $z=p_M(z')$. Since $D_M(\gamma^l\gamma'(z'))=D_M(w'), w'\in V'$ and (i) is satisfied, there exists $\alpha\in\pi_1(M)$ such $\alpha(z')=w'$. This implies that every neighborhood $V$ of $x$ contains at least two distinct plaques of ${\cal F}_M^U$ which are subsets of $\cup_{i=1,...,p}F_i$. This is impossible since each $F_i$ is compact. (See Haefliger [18] proposition 3.1).

\medskip

{\bf Proposition 4.7.}
{\it Under the hypothesis and notations of the proposition 4.6, a leaf $F$ is compact if and only if there exists $i\in\{1,...,p\}$ such that $F=F_i$ and   the holonomy of every non compact leaf is trivial.}

\medskip

{\bf Proof.}
 Consider a leaf $F$ distinct of $F_1,...,F_p$. Let $\hat F$ be a leaf of ${\cal F}_{\hat M}^U$ such that $p_M(\hat F)=F$; $D_M(\hat F)$ is contained in an affine subspace $w+U$, write $w=w_ne_n,w_n\neq 0$. Let $\gamma\in\pi_1(F)$, since $h_M(\gamma)(w+U)=w+U$, we deduce that $h_M(\gamma)(w)\in w+U$. We have: $h_M(\gamma)(w)=L(h_M)(\gamma)(w_ne_n)+h_M(\gamma)(0)=w+u',u'\in U$. The proposition 4.6 implies that we can suppose that $h_M(\gamma)(0)\in U$, we deduce that $L(h_M)(\gamma)(e_n)=e_n+v,v\in U$. This implies that $F$ does not have holonomy. If $F$ is compact, this  contradicts  the proposition 4.5.

\medskip

{\bf Remarks 4.3.}
The quotient of $\mathbb{R}^2-\{0\}$ by $h(x)=2x$ is the torus $T^2$. The linear foliation of $\mathbb{R}^2$ whose leaves are parallel to the $x$-axis  has two compact leaves which correspond to the images of $\{(x,0)\in\mathbb{R}-\{0\},x>0\}$ and $\{(x,0)\in\mathbb{R}-\{0\},x<0\}$ by the  covering map.

Sullivan and Thurston  [29] have endowed the three dimensional torus with the structure of an affine manifold such that the image of the developing map is $\mathbb{R}^3-\{D_1,D_2,D_3\}$ where $D_1,D_2,D_3$ are $3$ distinct lines which meet at the origin. This implies that the developing map of this affine structure is not a covering map.

\medskip

{\bf Proposition 4.8.}
{\it Let $(M,\nabla_M)$ be a compact,  complete affine manifold endowed with a linear foliation ${\cal F}_M^U$ of codimension $1$. Suppose that ${\cal F}_M^U$ has a compact leaf.

1. Then  all the leaves of ${\cal F}_M^U$ are compact.

2. There exists a finite cover $(M',\nabla_{M'})$ of $(M,\nabla_M)$ and an affine map $f:(M',\nabla_{M'})\rightarrow S^1$, here  $S^1$ is the quotient of $\mathbb{R}$ by a translation.}

\medskip

{\bf Proof.}
Let $ F_0$ be a compact leaf of ${\cal F}_M^U$. It is enough to show that the holonomy of every element of $\pi_1({ F}_0)$ is finite, and use  proposition 4.5 which implies that there exists  a finite cover $(M',\nabla_{M'})$ of $(M,\nabla_M)$ and an affine map $f:(M',\nabla_{M'})\rightarrow S^1$ whose fibres are leaves of the pullback of ${\cal F}_M^U$ to $M'$. Remark that since $(M',\nabla_{M'})$ is complete, the affine structure of the space of leaves $S^1$ is complete. We deduce that $S^1$ is the quotient of $\mathbb{R}$ by a translation.

Suppose that the holonomy of $F_0$ is not finite, without restricting the generality, we can suppose that $p_M(U)=F_0$. The proposition 4.6 implies that we can suppose that for every element $\gamma$ of $\pi_1(M)$, $h_M(\gamma)(0)$ is an element of  $U$. 
This contradicts the fact that the affine holonomy  representation of a compact and complete affine manifold is irreducible. (see Fried, Goldman and Hirsch [13] p. 496.)

\medskip

{\bf Theorem 4.1.}
{\it The global holonomy of a transversely oriented linear foliation  ${\cal F}_M^U$, of codimension strictly superior to $0$, defined on a compact affine manifold $(M,\nabla_M)$ is not trivial. }

\medskip

{\bf Proof.}
Suppose that the global holonomy $h_{{\cal F}^U_M}:\pi_1(M)\rightarrow Aff(\mathbb{R}^n/U)$ is trivial. Let $h_U:\mathbb{R}^n\rightarrow \mathbb{R}^n/U$ be the quotient map. There exists a map $g:M\rightarrow \mathbb{R}^n/U$ such that for every $x\in M, x'\in \hat M$ such that $p_M(x')=x, g(x)=h_U(D_M(x'))$. Since $h_U$ and $D_M$ are open maps, we deduce that $g$ is an open map. This implies that the image of $g$ is open. The image of $g$ is also closed since $M$ is compact, we deduce that the image of $g$ is $\mathbb{R}^n/U$. Contradiction.
  
\medskip

{\bf Corollary 4.3.}
{\it Let $(M,\nabla_M)$ be a compact affine manifold endowed with a linear foliation ${\cal F}_M^U$ of codimension $1$.  Then there exists a finite cover $M'$ of $M$ such that $H^1(M',\mathbb{R})\neq 0$.}

\medskip

{\bf Proof.}
Let $\gamma\in\pi_1(M)$, write $L(h_M)(\gamma)(e_n)=a_{\gamma}e_n+u_{\gamma}$ where $u_{\gamma}\in U$. The map defined on $\pi_1(M)$ by $a(\gamma)=a_{\gamma}$ is a morphism of groups. Up to a finite cover, we can suppose that for every $\gamma\in\pi_1(M)$, $a_{\gamma}>0$. We deduce that $Log(a_{\gamma})$ is a $1$-cocycle. Suppose that the cocycle  $Log(a_{\gamma})$ is trivial,   for every $\gamma\in \pi_1(M)$, $L(h_M)(\gamma)(e_n)=e_n+u_{\gamma}, u_{\gamma}\in U$. Write $h_M(\gamma)(0)=(b^{\gamma}_1,...,b^{\gamma}_n)$ in the basis $(e_1,...,e_n)$. The correspondence defined on $\pi_1(M)$ by $b(\gamma)= b^{\gamma}_n$ is a $1$-cocycle.  If this $1$-cocycle is trivial,  for every element $\gamma\in\pi_1(M)$, $h_M(\gamma)(0)=(b^{\gamma}_1,...,b^{\gamma}_{n-1},0)$. The theorem 4.1 shows that this is not possible.

\medskip

{\bf Remark 4.4.}
Let $\omega$ be a closed $1$-form defined on the manifold $M$. We can associate to $M$ a $1$-Cech cocycle defined as follows: consider a good  cover  $(U_i)_{i\in I}$ of $M$. Recall that for every finite subset $\{i_1,...,i_p\}\subset I$, the intersection $U_{i_1}\cap...\cap U_{i_p}$ is contractible. Let $\omega_i$ be the restriction of $\omega$ to $U_i$, there exists a differentiable function $f_i$ defined on $U_i$ such that $df_i=\omega_i$; we write $a_{ij}=f_j-f_i$. It is a constant function and the family of $(a_{ij})_{i,j\in I}$ defines a $1$-Cech cocyle of $\mathbb{R}_M$, the sheaf of locally constant real valued functions defined on $M$. This correspondence, which associates to $\omega$ the Cech cocycle $(a_{ij})_{i,j\in I}$ induces an isomorphism between the De Rham cohomology group $H^1_{DR}(M)$ and   $H^1_{Cech}(M,\mathbb{R}_M)$, the cohomology group of the sheaf $\mathbb{R}_M$.

 To an element $[a]$ of $H^1_{Cech}(M,\mathbb{R}_M)$ represented by the cocycle $(a_{ij})$ in the good covering $(V_i)_{i\in I}$, we can associate a  representation $h_a$ of $\pi_1(M)$ as follows: 
 
 Let $x$ be an element of $M$ and $\gamma$ an element of $\pi_1(M)$. Consider a path $c:[0,1]\rightarrow M$ whose homotopy class is $\gamma$. There exists a partition $[t_0=0,t_1],..,[t_k,t_{k+1}],[t_{k^{\gamma}-1},1]$ of $[0,1]$ such that the
   image of the restriction $c_{\mid[t_k,t_{k+1}]}$ is contained in $V_{i_k}$. We fix $V_{i_0}$ and suppose that for every $\gamma$, $c_{\mid [t_0=0,t_1]}\subset V_{i_0}$.   
   We  associate to $(a_{ij})_{i,j\in I}$ the morphism defined  by $h_a(\gamma)=a_{i_0i_1}+..+a_{i_{k^{\gamma}-1}i_1}$. The map which assigns $h_a$ to $(a_{ij})_{i,j\in I}$  induces an isomorphism between $H^1_{Cech}(M,\mathbb{R}_M)$ and $Hom(\pi_1(M),\mathbb{R})$.
 
 \bigskip
 
 Let $M$ be a manifold endowed with a transversely affine foliation ${\cal F}$. There exists an atlas $(U_i)_{i\in I}$, a differentiable map $k_i:U_i\rightarrow \mathbb{R}$ such that on $U_i\cap U_j$, $k_i=a_{ij}k_j+b_{ij}$ where $a_{ij}$ and $b_{ij}$ are constant. The restriction of ${\cal F}$ to $U_i$ is defined by the equation $dk_i=0$. The foliation ${\cal F}$ is also defined globally by a $1$-form $\omega$ without singularity. Let $\omega_i$ be the restriction of $\omega$ to $U_i$. There exists a differentiable function $f_i$ such that: $\omega_i={{dk_i}\over{f_i}}$.
 
 \bigskip
 
 We recall now the cohomology class associated to a transversely affine foliation defined by Bobo Seke in [5]. p. 7.
 
 \medskip
 {\bf Definition 4.1.}
 {\it  The form $\omega_1$ whose restriction to $U_i$ is ${{df_i}\over{f_i}}$ is well-defined and closed. Its cohomology class is called the Bobo Seke class of the transversely affine foliation ${\cal F}$.}
 
 \medskip
 We can show the following result:
 
\medskip

{\bf Proposition 4.9.}
{\it The cocycle $Log(a_{\gamma})$ defined  on the finite cover $M'$ of $M$ in the proof of the corollary 4.3 corresponds to the Bobo Seke class of ${\cal F}_{M'}^U$.}

\medskip

{\bf Proof.}
We suppose that ${\cal F}_M^U$ is transversally oriented, this enables to suppose that $M=M'$. Let $(V_i,\phi_i)_{i\in I}$ be an affine atlas such that the restriction of ${\cal F}_M^U$ to  $V_i$, is defined by $dx^i_n=0$ in the affine coordinates $(x_1^i,..,x_n^i)$. We have $dx^j_n=a_{ij}dx^i_n$ on $V_i\cap V_j$. The remark 4.4 shows  that the Cech cocycle $Log(a_{ij})$ is equivalent to  the cocycle $Log(a_{\gamma})$ which appears in the proof of the corollary 4.3. The distribution tangent to the foliation ${\cal F}_M^U$ is defined globally by a $1$-form $\omega$, we denote by $\omega_i$ the restriction of $\omega$ to $V_i$, there exists a differentiable function $f_i$ defined on $V_i$ such that $\omega_i={{dx_n^i}\over{f_i}}$. The fact that $dx^j_n=a_{ij}dx^i_n$ implies that $f_j=a_{ij}f_i$. We deduce that $Log(f_j)-Log(f_i)=Log(a_{ij})$. Since the Bobo Seke class is the cohomology class of the $1$-form $\omega_1$ whose restriction to $U_i$ is  ${{df_i}\over{f_i}}$, we  deduce that
$Log(f_j)-Log(f_i)=Log(a_{ij})$ is a Cech cocycle which represents  $\omega_1$.

\bigskip

{\bf 5. Codimension $1$ linear foliations defined on $3$-dimensional affine manifolds.}

\medskip

In this part, we are going to study codimension $1$ linear foliations defined on $3$-dimensional affine manifolds.

\medskip

{\bf Proposition 5.1.}
{\it Let $(M,\nabla_M)$ be a $3$-dimensional compact affine manifold endowed with a transversely, oriented, codimension $1$ linear foliation ${\cal F}_M^U$. Suppose that the developing map $D_M$ is injective, and ${\cal F}_M^U$  does not have a compact leaf. We suppose also that there exist a leaf $F_0$ of ${\cal F}_M^U$, a leaf $\hat F_0$ of ${\cal F}_{\hat M}^U$ above $F_0$, an element $x_0\in \mathbb{R}^3$ such that $\hat F_0=(x_0+U)\cap D_M(\hat M)$.  Then a finite cover of  $M$ is the total space of a bundle over $S^1$.}

\medskip

Firstly, we prove the following lemma:

\medskip

{\bf Lemma 5.1.}
{\it Let $G$ be a Lie subgroup of $Aff(\mathbb{R}^2)$, the group of affine transformations of $\mathbb{R}^2$. Suppose that $L(G)$, the linear part of $G$ contains $Sl(2,\mathbb{R})$. Then $G$ acts transitively on $\mathbb{R}^2$ or has a unique fixed point.}

\medskip

{\bf Proof.}
Let $L:G\rightarrow Gl(\mathbb{R}^2)$ defined by $L(A,a)=A$. Let $G_0$ the connected component of $L^{-1}(Sl(2,\mathbb{R}))$ and ${\cal G}_0$ its Lie algebra. Let $l:{\cal G}_0\rightarrow \mathbb{R}^2$ defined by $l(A,a)=a$, if $l$ is injective, it is a $1$-cocycle. The Whitehead's lemma implies that the cohomology class of $l$ is trivial and $G_0$ has a fixed point $x_0$ that we assume to be the origin, we deduce that $G$ acts transitively on $\mathbb{R}^2-\{0\}$ since $SL(2,\mathbb{R})$ acts transitively on $\mathbb{R}^2-\{0\}$. If $l$ is not injective, there exists $(0,a)\in {\cal G}_0, a\neq 0$. Let $(B,b)$ be another element of ${\cal G}_0$, $[(B,b);(0,a)]=(0,B(a))$ there exists $B\in L({\cal G}_0)=sl(2,\mathbb{R})$ such that $(a,B(a))$ is a basis of $\mathbb{R}^2$, we deduce that $G_0$ contains $(I,ta+t'B(a)), t,t'\in\mathbb{R}$ and acts transitively on $\mathbb{R}^2$.

\medskip

{\bf Proof of the proposition 5.1.}
The manifold $M$ is the quotient of a simply connected open subset $V$ of $\mathbb{R}^3$ by $\pi_1(M)$ since $D_M$ is injective. The theorem 9.1 of [24]   implies that $V$ is contractible since we have supposed that ${\cal F}_M^U$ does not have a compact leaf. The first paragraph p.26 [24] implies that the leaves of ${\cal F}_{\hat M}^U$ are contractible.

Suppose that $\pi_1(M)$ is not solvable. We have an exact sequence: 
$$
1\rightarrow Ker(h_{{\cal F}_M^U})\rightarrow\pi_1(M)\rightarrow Im(h_{{\cal F}_M^U})\rightarrow 1.\leqno(1)
$$ 

Since $Im(h_{{\cal F}_M^U})$ is a subgroup of the solvable group $Aff(\mathbb{R})$, it is solvable. We deduce that   $Ker(h_{{\cal F}_M^U})$ is not solvable. Write $U_F=x_0+U$ Consider the restriction map $r:Ker(h_{{\cal F}_M^U})\rightarrow Aff(U_F)$. Remark that the kernel of $r$ is an unipotent group, we deduce that the image $H$ of $r$ is not solvable.

Let $cd(H)$ be the cohomological dimension of $H$.

$\hat F_0$ is stable by $H$, since $\hat F_0=(x_0+U)\cap D_M(\hat M)$.
Remark that $cd(H)$  cannot be superior to $3$, since $H$ acts properly and freely on the contractible manifold $\hat F_0$, thus  the quotient space $\hat F_0/H$ is a $2$-dimensional  $K(H,1)$ Eilenberg McLane space.  If $cd(H)=2$, then $\hat F_0/H$ which is a covering of $F_0$ is compact. We deduce that $F_0$ is compact. Contradiction. We deduce that $cd(H)=1$. The theorem of Stallings (Serre [27] example c) p. 90) implies that $H$ is a free group. 
  The Zariski closure $Z(h_{F_0})$, of the image of $h_{F_0}$, is a subgroup of $Aff(U_F)$ whose linear part contains $Sl(2,\mathbb{R})$, since it contains  $Ker(h_{{\cal F}_M^U})$ which is not solvable.  Remark that $Z(h_{F_0})$ preserves  the boundary of  $V\cap \hat F_0$, since it is an algebraic set, (see Sullivan Thurston [29] p.22). The lemma 5.1 implies that  $V\cap \hat F_0$ is isomorphic as an affine space to  $\mathbb{R}^2$ or  $\mathbb{R}^2-\{0\}$. It cannot be isomorphic to $\mathbb{R}^2-\{0\}$ since $\hat F_0$ is contractible.  We deduce that $F_0$ is endowed with the structure of a complete affine manifold. This implies that the action of $h_{F_0}$ is free and every element in the image of $L(h_{F_0})$ has $1$ as an eigenvalue. This implies that the linear part of every element of $Z(h_{F_0})$ has $1$ as an eigenvalue, contradiction, since $L(Z(h_{F}))$ contains $Sl(2,\mathbb{R})$.

We deduce that $\pi_1(M)$ is solvable. The lemma 3.1 implies that a finite cover of $M$ is the total space of a bundle over $S^1$.
\medskip

 The propositions 5.1  enables to show:

\medskip

{\bf Theorem 5.1.}
{\it Let $(M,\nabla_M)$ be a $3$-dimensional compact affine manifold, whose developing map is injective. Suppose that $(M,\nabla_M)$ is endowed with a codimension $1$ linear foliation ${\cal F}_M^U$, such that there exist a leaf $F_0$ of ${\cal F}_M^U$, a leaf $\hat F_0$ of ${\cal F}_{\hat M}^U$ above $F_0$, an element $x_0\in \mathbb{R}^3$, such that $\hat F_0=(x_0+U)\cap D_M(\hat M)$. Moreover, suppose that $F_0$ is compact or all the leaves of ${\cal F}_M^U$ are not compact,  then a finite cover of  $M$ is the total space of a bundle over $S^1$.}

\medskip

{\bf Proof.}
Under the hypothesis of the theorem 5.1, suppose that $F_0$ is compact, the theorem of Benzecri implies that $\pi_1(F_0)$ is polycyclic. The fact that $\hat F_0=(x_0+U)\cap D_M(\hat M)$ implies that $\hat F_0$ is stable by $ker(h_{{\cal F}_M^U})$. By applying the proposition 4.1 and the lemma 3.1, we deduce that a finite cover of $M$ is the total space of a bundle over $S^1$. If all the leaves of ${\cal F}_M^U$ are not compact, we apply proposition 5.1.

\medskip

{\bf Corollary 5.1}
{\it Let $(M,\nabla_M)$ be a $3$-dimensional affine manifold endowed with a codimension $1$ linear foliation. Suppose that the developing map of $M$ is injective and  its image is convex. Then a finite cover of $M$ is a total space of a bundle over $S^1$.}

\medskip

{\bf Proof.}
Let $F_0$ be a leaf of ${\cal F}_M^U$. Consider a leaf $\hat F_0$ above $F_0$, and $x_0$ an element of $\hat F_0$. Since $D_M(\hat M)$,  is convex, we deduce that $\hat F_0= (x_0+U)\cap D_M(\hat M)$. We can apply the theorem 5.1 to deduce that a finite cover of $M$ is the total space of a bundle over $S^1$.

\medskip

{\bf Remark 5.1.}

In the proof of the proposition 5.1, we have used the fact that the fundamental group of a complete $2$-dimensional affine manifold is solvable. There exist $2$-dimensional affine manifolds whose fundamental group is not solvable: for example,  $\mathbb{R}^2-\{(0,0); (0,1)\}$ is a $2$-dimensional affine manifold whose fundamental group is the free group.

  Fried and Goldman [12] have classified $3$-dimensional geodesically complete closed affine manifolds. They have shown that their fundamental group contains a solvable subgroup of finite index. More generally,
Auslander has conjectured that the fundamental group of a closed geodesically complete affine manifold contains a polycyclic subgroup of finite index. In, [23] Milnor has shown that a group without torsion, and which contains a polycyclic group of finite index is the fundamental group of a complete affine manifold, and he has asked whether the fundamental group of every complete manifold is virtually polycyclic. Let $F_2$ be the free group generated by $2$ elements. In [22] Margulis has constructed an example of properly discontinuous action of $F_2$ by affine transformations  on $\mathbb{R}^3$. The quotient of $\mathbb{R}^3$ by this action of $F_2$ is an handlebody. Drumm [9] has constructed fundamental domains for this action by defining the notion of crooked plane see also [7].

\bigskip

 \medskip
 
 {\bf Acknowledgements.}
 
 The author thanks the referee whose remarks have considerably improved the presentation of this paper. 

\medskip

This research did not receive any specific grant from funding agencies in the public, commercial, or
not-for-profit sectors.

\bigskip

{\bf References.}

\medskip

[1] A. Albert, Structure of algebras, Vol. 24. American Math. Soc. 1939.

\smallskip

[2] M. Aschenbrenner, S. Friedl, H. Wilton,  3 manifold groups. https://arxiv.org/pdf/1205.0202.pdf

\smallskip

[3]J.P. Benz\'ecri,  Sur les variétés localement affines et localement projectives,  Bulletin de la Société Mathématique de France, 88 (1960) 229-332.

\smallskip

[4] R. Bieri,   Homological dimension of discrete groups, Queen Mary College Mathematics Notes, 1976.

\smallskip

[5] Bobo Seke, Sur les feuiffetages transversalement affines de codimension $1$,
 Annales de l'Institut Fourier, 30 (1980) 1-29.

\smallskip

[6] Y. Carri\'ere, Yves, Autour de la conjecture de L. Markus sur les variétés affines, Inventiones mathematicae, 95  (1989): 615-628.

\smallskip

[7] V. Charette, T.  Drumm, W.  Goldman W, Affine deformations of a three-holed sphere,  Geometry and Topology, 14 (2010) 1355-1382.

\smallskip

[8] S. Choi, The decomposition and classification of radiant affine 3-manifolds 154, American Math. Soc., 2001.

\smallskip

[9] T. Drumm,  Fundamental polyhedra for Margulis space-times,  Topology,  (1992) 677-683.

\smallskip

 [10] F.T.  Farrell, L.E.  Jones,  Rigidity for aspherical manifolds with $\pi_1\subset GL_m (\mathbb {R}) $, Asian Journal of Mathematics, 2 (1998): 215-262.

\smallskip

[11] M. H. Freedman,  F. Topology of 4-Manifolds  Volume 39, Princeton University Press, 2014

\smallskip

 [12] D. Fried, W.   Goldman,  Three-dimensional affine crystallographic groups,  Advances in Mathematics, 47 (1983) 1-49.

\smallskip

[13] D. Fried, W. Goldman, M.  Hirsch, M, Affine manifolds with nilpotent holonomy,  Commentarii Mathematici Helvetici, 56 (1981) 487-523.

\smallskip

[14] E. Ghys, V.  Sergiescu,  Stabilité et conjugaison différentiable pour certains feuilletages, Topology, 19  (1980) 179-197.

\smallskip

[15] W. Goldman, M.  Hirsch,  The radiance obstruction and parallel forms on affine manifolds, Transactions of the American Mathematical Society, 286 (1984) 629-649.

\smallskip

[16] W. Goldman, M. Hirsch,  Affine manifolds and orbits of algebraic groups, Transactions of the American Mathematical Society, 295 (1986) 175-198.

\smallskip

[17] A. Haefliger,  Structures feuillet\'ees et cohomologie \`a valeur dans un faisceau de groupoides, Commentarii Mathematici Helvetici, 32 (1958) 248-329.

\smallskip

[18] A. Haefliger,  Vari\'et\'es feuillet\'ees, Annali della Scuola Normale Superiore di Pisa Classe di Scienze,   16 (1962) p. 367-397.

\smallskip

[19] G. Hector, H.  Hirsch,  Introduction to the Geometry of Foliations: Part B, F. Vieweg  John, 1983.

\smallskip

[20] W.C. Hsiang, C.T.C.  Wall,  On homotopy tori II,  Bull. London Math. Soc.,  (1969) 341-342.

\smallskip

[21] L.  Markus,  Cosmological models in differential geometry, U. Minesota, 1962.

\smallskip

[22] G. Margulis,    Free  properly  discontinuous  groups  of  affine  transformations,  Dokl. Akad. Nauk SSSR272, (1983) 937–940.

\smallskip

[23] J. Milnor,  On fundamental groups of complete affinely flat manifolds, Advances in Mathematics, 25 (1977) 178-187.

\smallskip

[24] S. Novikov,  The topology of foliations, Trudy Moskov. Mat. Obshch, 1965  248–278. http://www.mi-ras.ru/~snovikov/23.pdf

\smallskip

[25] M.S. Raghunathan,  Discrete subgroups of Lie groups, Springer, 1972.

\smallskip

[26] H. Rosenberg,   Foliations by planes,  Topology, 7 (1968) 131-138.

\smallskip

[27] J.P. Serre,  Cohomologie des groupes discrets,  Prospects in mathematics, (Proceedings of Symposium, Princeton University, Princeton, NJ, 1970), (1971) 77-169.

\smallskip

[28] J. Smillie,  An obstruction to the existence of affine structures,  Inventiones mathematicae, 64 (1981) 411-415.

\smallskip

[29] D. Sullivan, W.  Thurston,  Manifolds with canonical coordinate charts: some examples, Enseign. Math., 29 (1983) 15-25.

\smallskip

[30] D. Tischler,  On fibering certain foliated manifolds over $S^1$,  Topology, 9 (1970)  153-154.

\smallskip

[31] A. Tsemo,  Dynamique des vari\'et\'es affines,   Journal of the London Mathematical Society, 63 (2001 469-486.

\end{document}